\documentclass[11pt,english]{article}

\usepackage{amsmath}

\usepackage{amsfonts}
\usepackage{graphicx}
\usepackage{amssymb}
\usepackage{caption}
\usepackage{subcaption}
\usepackage[english]{babel}
\usepackage[thmmarks]{ntheorem}

\usepackage{gensymb}
\usepackage[left=1cm,right=1cm,top=1cm,bottom=2cm]{geometry}
\usepackage{wrapfig}
\usepackage{mathenv}
\usepackage[T1]{fontenc}
\usepackage{lmodern}
\usepackage{dsfont}
\usepackage{enumitem}
\usepackage{mathtools}
\usepackage{stmaryrd}
\usepackage{hyperref}

\theorembodyfont{\normalfont}
\newtheorem{rmq}{Remark}[section]
\theoremseparator{:}

\theoremseparator{:}

\theoremstyle{plain}
\theoremheaderfont{\normalfont\bfseries}\theorembodyfont{\slshape}
\theoremseparator{:}
\newtheorem{theorem}{Theorem}

\theoremstyle{plain}
\theoremheaderfont{\normalfont\bfseries}\theorembodyfont{\slshape}
\theoremseparator{:}

\theoremstyle{plain}
\newtheorem{lemme}{Lemma}[section]

\newtheorem{prop}{Proposition}[section]

\theoremstyle{plain}
\theoremseparator{.}

\theoremnumbering{arabic}
\newtheorem{cor}{Corollary}[section]

\theoremstyle{nonumberplain}
\theoremheaderfont{\itshape}
\theorembodyfont{\normalfont}
\theoremsymbol{\ensuremath\square}

\theoremstyle{empty}
\theoremheaderfont{\itshape}
\theorembodyfont{\normalfont}
\theoremsymbol{\ensuremath\square}
\newtheorem{proof2}{}

\usepackage{titlesec}
\titlespacing{\section}{15pt}{*4}{*1.5}
\titlespacing{\subsection}{30pt}{*4}{*1.5}
\titlespacing{\subsubsection}{45pt}{*4}{*1.5}

\newcommand{\itg}{\displaystyle\int }

\newcommand{\R}{\mathbb{R}}
\newcommand{\C}{\mathbb{C}}
\newcommand{\T}{\mathcal{T}}
\newcommand{\cont}{\mathcal{C}}
\newcommand{\G}{\mathcal{G}}
\newcommand{\pinf}{^{\infty}}

\newcommand{\s}{\mathcal{S}}
\newcommand{\Spd}{\mathbb{S}^{d-1}}

\newcommand{\dmu}{_{d-1}}
\newcommand{\dpu}{_{d+1}}
\newcommand{\M}{\mathcal{M}}

\DeclareMathOperator{\supp}{Supp}

\author{\textsc{ Corentin Gentil}\footnote{École Normale Supérieure, 45 Rue d'Ulm, corentin.gentil@ens.fr},  \textsc{ Côme Tabary}\footnote{École Normale Supérieure, 45 Rue d'Ulm, come.tabary@ens.fr} }
\date{Juillet 2021}
\title{\text{A remark on restriction theorems and their application to Strichartz estimates}}

\begin{document}

\maketitle
\section*{Abstract}
We present an elementary approach to prove restriction theorems for particular surfaces for which the Tomas-Stein theorem does not apply, which in turn provide short proofs for well-known Strichartz estimates for associated PDEs. The method consists in applying simple restriction theorems to the level sets of the surfaces studied (for instance, spheres in the case of the cone) and then integrating the inequalities over all level sets. This allows for a different proof of sharp estimates for the wave equation, and of some estimates related to the Euler equations in the rotational framework. 

\section*{Résumé}
Nous présentons une approche élémentaire pour prouver des théorèmes de restriction sur certaines surfaces pour lesquelles le théorème de Tomas-Stein ne s'applique pas, qui permettent de retrouver des inégalités de Strichartz déjà connues. Notre méthode consiste à appliquer des théorèmes de restriction simples aux courbes de niveaux d'une surface donnée (par exemple, des sphères dans le cas du cône), puis à intégrer ces inégalités. Cela permet une preuve plus courte des estimées de Strichartz pour l'équation des ondes, et de certaines estimées pour l'équation d'Euler en rotation rapide.
\\
\\
\\
\begin{center}
\large{\textsc{Acknowledgments}}
\end{center}

\begin{center}
	\begin{minipage}{16cm}
		We are very grateful to Isabelle Gallagher for all her advice during the writing of this note and for her ideas that gave birth to this work. We also wish to thank the Département de Mathématiques et Applications of the ENS for their warm welcome.
	\end{minipage}
\end{center}
\medskip
\section{Introduction}
\subsection{Restriction theory}

For any function $f$ in the Lebesgue space $L^1(\R^d)$, its Fourier transform $\hat{f}$, defined for $\xi \in \R^d$ as
$$\hat f (\xi) = \int_{\R^d}e^{-ix\cdot\xi}f(x)dx,$$
is continuous, so the restriction of $\hat f$ to any surface of $\R^d$ is well defined. If $f$ is not integrable but only lies in some $L^p$ space, its Fourier transform can still be defined and the Hausdorff-Young inequality implies that if $1\leq p \leq 2$, then $\hat f$ is in the dual space $L^{p'}$, where $\frac{1}{p}+\frac{1}{p'}=1$. Thus $\hat f$ is only defined almost everywhere, and it is not clear whether it is possible to meaningfully restrict $\hat{f}$ to a hypersurface $S$ of $\R^d$, since $S$ is a set of measure zero. The restriction problem asks for which exponents $p$ and $q$ a restriction estimate of the form
$$ \Vert \hat{f}\vert_{S}\Vert_{L^q(S,d\sigma)} \leq C \Vert f \Vert_{L^p(\R^d)}$$ holds for all functions $f$ in the Schwartz class $\s(\R^d)$, where $d\sigma$ is a measure on $S$. Such estimates indeed allow one to define the restriction of the Fourier transform $f\mapsto \hat{f}\vert_{S}$ as a bounded operator from $L^p(\R^d)$ to $L^q(S,d\sigma)$.\\
Tomas and Stein discovered \cite{tsb}, \cite{ts} that this is possible for surfaces which are "sufficiently curved", a result that is now known as the Tomas-Stein theorem:

\begin{theorem}[\cite{tsb}, \cite{ts}] 
Let $d\geq 2$ and $p_0 = \frac{2(d+1)}{d+3}$. Let $S$ be a smooth and compact hypersurface in $\R^d$ endowed with a smooth measure $d\sigma$. Suppose $S$ has non vanishing Gaussian curvature at every point. Then for all $p\in [1, p_0]$, there exists a positive constant $C$ such that for any $f\in\s(\R^d)$:
$$ \Vert \hat f \vert_{S} \Vert_{L^2(S,d\sigma)} \leq C \Vert f \Vert_{L^p(\R^d)}.$$
\end{theorem}
The canonical example of such a surface is the sphere $\Spd$ of $\R^d$. When the hypothesis of curvature or of compactness is removed, obtaining restriction estimates does not become impossible, but the range of possible exponents usually gets smaller, as we will illustrate with the surfaces considered in this note. A first step towards understanding the necessity of curvature is to consider a simple counter-example (found for instance in \cite{tao}) preventing any non-trivial restriction estimate on a hyperplane (endowed with the natural surface measure): consider $\varphi\in \s (\R^{d-1})$, and define $h:\R^d\mapsto\C$ by:
$$h(x)=\dfrac{\varphi(x')}{\sqrt{1+x_{d}^{2}}}, \; \; \; x=(x',x_{d}).$$
Then $h$ belongs to every $L^p$ space for $p>1$, but its Fourier transform is not well defined on the hyperplane ${\xi_d=0}$, or even on any bounded subset (with positive surface measure) of this hyperplane. Curvature is thus mandatory to obtain restriction estimates. Lack of compactness, on the other hand, can sometimes be dealt with using scaling arguments, provided the surface possesses some symmetry.
\\
Restriction theory has many applications to other topics, from number theory to the Bochner-Riesz and Kakeya conjectures \cite{tao}. It can also be used in the study of partial differential equations, especially to provide Strichartz estimates. This is the focus of the following section.

\subsection{Strichartz estimates}

Considering an evolution partial differential equation, a Strichartz estimate is a bound for the space-time norm of a solution $u$ by the norm of the initial datum $u_0$, that is to say an inequality of the type
$$ \Vert u(t,x) \Vert_{L^p_tL^q_x} \leq C \Vert u_0 \Vert_{L^r_x}.$$
Here, $p, q, r$ are suitable exponents and $ \Vert u \Vert_{L^p_t L^q_x  }$ is the standard Lebesgue space-time norm on $L^p(\R,\ L^q(\R^d))$.\\
Strichartz inequalities are commonly obtained using a dispersive estimate. It is for instance the case for the linear Schrödinger equation:
$${(S):{\left\{\begin{array}{lrcl}
i\partial_t u + \Delta u =0, \ \ \ \text{for   } (t,x)\in\R^+ \times \R^d \\
u\vert_{t=0} = u_0. \end{array}\right.}}$$
The solution $u$ is, for $t>0$, 
$$ u(t,\cdot) = \frac{1}{(4\pi it)^{\frac d 2}} e^{\frac{\vert x \vert^2}{4t} } \star u_0.$$
This formula leads to the following dispersive estimate:
$$ \Vert u(t,\cdot) \Vert_{L\pinf_x(\R^d)} \leq C \frac{1}{(4\pi t)^{\frac d 2}} \Vert u_0 \Vert_{L^{1}_x(\R^d)},$$
for $ t >0$.\\
We can combine this inequality with a result from functional analysis known as the $TT^*$ argument (see \cite{velo}, \cite{kato}, \cite{keel}), to obtain Strichartz estimates for the Schrodinger equation (as it is done in \cite{bcg}): 
$$\Vert u \Vert_{L^q_t L^r_x} \leq C \Vert u_0 \Vert_{L^2_x},$$
with $q\in(2,\infty]$ and $r$ satisfying the scaling condition $ \frac 2 q + \frac d r = \frac d 2 $.\\
\\
However, once we have written the explicit formula of the solution, getting dispersive estimates is not necessarily immediate and it may require more advanced tools. Concerning the wave equation 
$$  (W):   {\left\{\begin{array}{lrcl}
\partial^2_{tt} u (t,x)- \Delta u(t,x) = 0, \ \ \ \ \ \ \ \ \ \text{for   } (t,x)\in\R^+ \times \R^d\\
u|_{t=0} = f\ \ \ \ \text{and}\ \ \ \ \partial_t u|_{t=0} = g, \end{array}\right.} $$
dispersive estimates are obtained using the stationary phase method, see for instance \cite{bcg} and \cite{cdgg}. It is well known that the solution can be written, once again using the Fourier transform:
$$ u(t,x) = \mathcal{F}^{-1} (e^{it\vert x\vert } \gamma_+(\xi) + e^{-it\vert x\vert } \gamma_-(\xi))$$
where $\gamma_{\pm} (\xi) =\frac{1}{2} \bigg(  \hat f(\xi) \pm \frac{\hat g(\xi)}{i\vert\xi\vert} \bigg)$. In the case $\gamma_{\pm}$ are supported in a ring, a stationary phase argument yields the following dispersive estimate:
$$ \Vert u(t,\cdot) \Vert_{L\pinf(\R^d)}\leq \frac{C}{\vert t \vert^{\frac{d-1}{2}}} (\Vert f \Vert_{L^1} + \Vert g \Vert_{L^1}).$$
This leads to the following Strichartz estimates for the wave equation:
$$ \Vert u \Vert_{L^q_t L^p_x} \leq C (\Vert \nabla f \Vert_{L^2} + \Vert g \Vert_{L^2})$$
with $p \in [2,\infty ],\ q\in [2,\infty )$ satisfying $ \frac{1}{q} + \frac{d}{p} = \frac{d}{2}-1$.\\
We now turn to the connection between Strichartz estimates and Fourier restriction theorems. As it was shown by Strichartz in \cite{stri}, solutions to some PDEs can be written as Fourier transforms on surfaces. Consider once again the solution of the free wave equation $(W)$ ; it can be written as:
$$ u(t,x) = \itg_{\R^d_{\xi} } e^{it\vert \xi\vert + i x \cdot \xi} \gamma_+(\xi) d\xi + \itg_{\R^d_{\xi}} e^{-it\vert \xi\vert + i x \cdot \xi} \gamma_-(\xi) d\xi.$$

The solution can thus be seen as the sum of two inverse Fourier transforms, each one on a cone $\mathcal{C}_{\pm} = \{ (\xi, \pm \vert\xi\vert ),\ \xi \in \R^d \}$. Indeed,

$$ u(t,x) = \itg_{\mathcal{C_+}} e^{i (t,x) \cdot \zeta } \ \ \Tilde{\gamma}_+(\zeta)  \ d\mu(\zeta) + \itg_{\mathcal{C_-}} e^{i (t,x) \cdot \zeta }\ \  \Tilde{\gamma}_-(\zeta)  \ d\mu(\zeta), $$
with $\Tilde{\gamma}_{\pm}(\xi, \pm \vert \xi \vert) = \gamma_{\pm}(\xi) $. The measure $d\mu(\xi, \pm \vert \xi \vert) = d\xi$ is simply the pullback of Lebesgue measure by the canonical projection $(\xi, \pm \vert\xi\vert) \mapsto \xi$. This is exactly the dual of the restriction operator $f\mapsto \hat{f}\vert_{\mathcal{C}_{\pm}}$ applied to $\Tilde{\gamma}_{\pm}$, thus inciting to resort to restriction theory.\\
Let us admit for a brief moment that we have the following restriction estimate on the cone: 

$$\left\Vert {\hat v}\vert_{\cont_{\pm}}\  {\vert \xi \vert^{-\frac{1}{2}}} \right\Vert_{L^2(\cont_{\pm},\ d\mu)} \leq C \Vert v \Vert_{L^{p_0}(\R^{d+1})}, $$

where $p_0 = \frac{2(d+1)}{d+3}$. In its dual form, it writes, for $w\in L^2(\cont_{\pm}, d\mu)$:

$$ \left\Vert \mathcal{F}^{-1} (\vert \xi \vert^{-\frac{1}{2}} w \ d\mu) \right\Vert_{L^{p_0'}(\R^{d+1})} \leq C \Vert w \Vert_{L^2(\cont_{\pm}, d\mu)}   .$$

In our case, setting $w = \Tilde{\gamma}_{\pm} \vert \xi \vert^{\frac{1}{2}}$, we notice $$ \left\Vert \Tilde{\gamma}_{\pm} \vert \xi \vert^{\frac{1}{2}} \right\Vert_{L^2(\cont_{\pm}, d\mu)} = \Vert \gamma_{\pm} \Vert_{\dot{H}^{\frac{1}{2}}(\R^d)}.  $$
Finally, applying the dual inequality and using the definition of $\gamma_{\pm}$, we obtain the Strichartz estimate
$$ \Vert u \Vert_{L^{p_0'}_{t,x}} \leq C \Vert f \Vert_{\dot{H}^{\frac{1}{2}}(\R^d)}  + \Vert g \Vert_{\dot{H}^{-\frac{1}{2}}(\R^d)} .$$

Hence we have seen that, through fairly simple rewriting of the representation formula, the obtention of the Strichartz estimate above reduces to the proof of the following restriction theorem:

\begin{theorem}[\cite{stri}]
Let $d\geq 2$ and consider the cone $\cont = \{ (\xi, \vert \xi \vert ),\ \xi \in \R^d \}\subset \R^{d+1}$ endowed with surface measure $d\sigma(\xi,|\xi|) = \frac{d\xi}{|\xi|}$. For $p_0= \frac{2(d+1)}{d+3}$, there exists a constant $C>0$ such that for any function $f\in \s(\R^{d+1})$:
$$ \Vert \hat{f}|_{\cont} \Vert_{L^2(\cont, d\sigma)} \leq C \Vert f \Vert_{L^{p_0}(\R^{d+1})}. $$
\end{theorem}

Note that Theorem 2 is not a direct consequence of Theorem 1 because the cone has zero Gaussian curvature. The original proof by Strichartz in \cite{stri}, for general quadratic surfaces, relies on complex analysis arguments. The first goal of this article is to give a different proof of this theorem. It consists in applying Theorem 1 to the spheres of increasing radius that make up the cone, an approach motivated by the fact that the exponent $p_0$ for the cone $\mathcal{C}\subset \R^{d+1}$ is the endpoint one for the sphere $\Spd$ of $\R^d$.\\

The resort to restriction theory to obtain Strichartz estimates is motivated by its efficiency even in settings where dispersion is not easily available. Such settings include bounded domains and compact manifolds, see for instance the work of Bourgain \cite{bourgain} for the Schrödinger equation on $\mathbb{T}^d$. This is also the case of the Heisenberg group $\mathbb{H}^d$, where there is no dispersion at all \cite{bgx}, but where Strichartz estimates for the wave equation and the Schrodinger equation were recently obtained through restriction theorems by Bahouri, Barilari and Gallagher \cite{bbg}.\\ 

In this article, we are also interested in using restriction theory in the study of the Euler equation in a rotational framework. The starting point is the following system:
$$  (EC_{\varepsilon}):   {\left\{\begin{array}{lrcl}
\partial_t u + u \cdot \nabla u + \dfrac{e_3 \wedge u}{\varepsilon} + \nabla p = 0 \\
\text{div} \ u = 0 \end{array}\right.} $$

with initial datum $u(0,\cdot)=u_0 $, and where the unknowns are both the vector field $u$, which is the velocity of the fluid, and the pressure field $p$. It is the zero-viscosity limit of the equations 
$$  (NSC_{\varepsilon}):   {\left\{\begin{array}{lrcl}
\partial_t u + u \cdot \nabla u - \nu \Delta u + \dfrac{e_3 \wedge u}{\varepsilon} + \nabla p = 0 \\
\text{div} \ u = 0, \end{array}\right.} $$
which are studied in details in Chapter 5 of \cite{cdgg}. To get rid of the unknown pressure field, one can apply the Leray projector $\mathbb{P}$ onto divergence-free vector fields, resulting in the following equations:
$$     {\left\{\begin{array}{lrcl}
\partial_t u +  \mathbb{P}\ (\dfrac{e_3 \wedge u}{\varepsilon}) = -\mathbb{P}\ (u \cdot \nabla u)  \\
\text{div} \ u = 0. \end{array}\right.} $$
We put the non-linear term $\mathbb{P}\ (u\cdot \nabla u)$ on the right-hand side because it is often treated as a source term (\cite{cdgg}, \cite{cor}) by first solving the homogeneous equation, and then using Duhamel formula:
$$ u(t) = \G^{\varepsilon} (t) u_0 - \itg_0^t \G^{\varepsilon} (t-t') \mathbb{P}\ (u(t') \cdot \nabla u(t')) dt' ,$$
where $\G^{\varepsilon}$ denotes the semigroup of the homogeneous equation
$$(C_{\varepsilon}): \partial_t u  + \mathbb{P}\ \dfrac{e_3 \wedge u}{\varepsilon} = 0.$$

Strichartz estimates for this semigroup were proven in \cite{dut} and \cite{cor}. They are in turn used to prove long-term existence of solutions to the Euler equations in the rotational framework. The proof in \cite{dut} adapts the arguments in \cite{cdgg} to the non-viscous case, using the Coriolis force to obtain a dispersive estimate ; whereas in \cite{cor} a slightly more precise study is made to obtain sharp results. \\
Our goal is to obtain such estimates through a global restriction theorem on a suitable surface. To do so, we first need to write the solution of $(C_{\varepsilon})$ as an inverse Fourier transform. We can reduce to the case $\varepsilon=1$ through a simple rescaling of $u$, and write the equation in Fourier space, which gives the linear ODE
$$ \partial_t \hat{u} (t,\xi) = A(\xi) \hat{u} (t, \xi),$$
where $A(\xi)$ is a matrix with eigenvalues $0$ and $ \pm i \frac{\xi_3}{\vert \xi \vert}$. Using that $u_0$ is divergence free, we infer that $\hat{u_0}$ has no component of eigenvalue $0$. By projecting $\hat{u_0}(\xi)$ onto the two other eigenspaces, the semigroup $\G^1$ can be written in terms of the semigroups $e^{\pm it\frac{D_3}{\vert D \vert}}$, defined as:
$$(e^{\pm it\frac{D_3}{\vert D \vert}} v)\ (x)= \itg_{\R^3}\hat{v}(\xi) e^{\pm i t \frac{\xi_3}{\vert \xi \vert} +i x\cdot \xi }d\xi.$$
Recognizing (up to a minus sign) an inverse Fourier transform on the surface
$$ \M:= \left\{ \bigg(\xi, \frac{\xi_3}{\vert \xi \vert}\bigg),\ \ \ \xi \in \R^3\right\} \subset \R^4,$$
Strichartz estimates for the semigroup $\G$ are available as soon as we have a restriction theorem on the surface $\M$, simply by using its dual form.\\

The second goal of this note is thus to prove a global restriction theorem on $\M$. Notice that this is once again not a direct application of Theorem 1, since $\M$ is not compact, and has vanishing Gaussian curvature. However, we can follow the same approach as for the cone. By introducing Sobolev weights in the surface measure to deal with non-compactness, and applying Theorem 2 to suitable sub-manifolds of $\M$ of lower dimension, we will prove the following theorem:
\begin{theorem}
Consider the surface $\M$ endowed with the measure $$d\sigma\bigg(\xi,\frac{\xi_3}{|\xi|}\bigg) = \frac{d\xi}{|\xi|^2}.$$ There exists a constant $C>0$ such that for any $f\in\s(\R^4)$, one has
$$ \Vert \hat f\vert_{\M}  \Vert_{L^2(\M,d\sigma)} \leq C \Vert f \Vert_{L^{\frac{6}{5}}(\R^4)}. $$
In terms of PDEs, this implies that for any $u_0 \in \s(\R^3),$ the following Strichartz estimate holds:
$$ \Vert e^{it\frac{D_3}{\vert D \vert}} u_0 \Vert_{L^6_{t} L^6_{x}} \leq C \Vert u_0 \Vert_{\dot H^1}. $$
\end{theorem}
The Strichartz inequalities obtained using this theorem are not sharp. More precisely, the optimal range of estimates obtained in \cite{cor} is
$$ \Vert e^{it\frac{D_3}{\vert D \vert}} u_0 \Vert_{L^q_{t}(\R, L^r_{x}(\R^d))} \leq C \Vert u_0 \Vert_{\dot H^{\frac{3}{2}-\frac{3}{r}}(\R^d)}, $$
for any $(q,r)\neq(2,\infty)$ such that $\frac{1}{q}+\frac{1}{r}\leq \frac{1}{2}$. Our technique only allows us to obtain, after interpolation with energy conservation, the range $\frac{1}{q}+\frac{1}{r}\leq \frac{1}{3}$, for $q\geq 6$. It is easy to identify where the sharpness is lost in the proof of Theorem 3, a point that we will further discuss in the next section and in Remark 3.1.\\
The Sobolev regularity required for the initial data in the Strichartz estimates can be explained by the fact that the surface looks like a hyperplane at infinity. We have seen that hyperplanes do not admit non-trivial restriction theorems, and it is somewhat similar in this case. It will be the object of Proposition 3.1.
\\
Let us conclude this section by noticing that Theorem 3 extends to a more general surface, namely
$$\mathcal{M}_{F}=\left\{ \bigg(\xi,\dfrac{\vert \xi \vert_{F}}{F\vert\xi\vert}\bigg), \xi \in \R^{3}\right\}\subset \R^4,$$
where $\vert \xi \vert_{F}^{2}=\vert\xi_{h}\vert^{2}+F^{2}\vert\xi_{3}\vert^{2}$ and $F\in(1,\infty]$  (note that $\M$ from Theorem 3 is simply the case $F=\infty$ up to an absolute value). This surface appears in \cite{charve} in the study of the primitive system, in which $F$ is the Froude number. The same estimates then hold for $\mathcal{M}_{F}$:

\begin{theorem}
Let $F\in (1,\infty).$ Consider the surface $\mathcal{M}_{F}$ with the measure $d\sigma(\xi,\frac{\xi_3}{|\xi|}) = \frac{d\xi}{|\xi|^2}$. There exists a constant $C>0$ such that for any $f\in\s(\R^4)$,
$$\Vert \hat f\vert_{\M}  \Vert_{L^2(\M,\frac{d\xi}{\vert\xi\vert^2})} \leq C \frac{F^{\frac 5 3}}{(F-1)^{\frac 5 9}} \Vert f \Vert_{L^{\frac{6}{5}}(\R^4)}. $$
\end{theorem}
This estimate could in turn provide new results on the primitive system with vanishing viscosity: we postpone this question to a future work.

\subsection{Motivation, layout and notation}

To introduce and motivate the method we will use in the proofs of Theorems 2, 3 and 4, we state the following lemma (mentioned in \cite{tao} as an exercise):

\begin{lemme}
Let $A \subset \R^m$ equipped with a $\sigma$-finite measure $d\sigma_A$ and $B\subset \R^n$ equipped with a $\sigma$-finite measure $d\sigma_B$. Assume there exist $p\leq q$ and $C>0$ such that for any Schwartz functions $f\in\s(\R^m)$ and $g\in\s(\R^n)$, $$ \Vert \hat{f}|_{A} \Vert_{L^q(d\sigma_A)} \leq C \Vert f \Vert_{L^p(\R^m)}\ \ \ \ \text{and} \ \ \ \ \Vert \hat{g}|_{B} \Vert_{L^q(d\sigma_B)} \leq C \Vert g \Vert_{L^p(\R^n)}.$$
Then for any Schwartz function $h\in \s(\R^{n+m})$, $$ \Vert \hat{h}|_{A\times B} \Vert_{L^q(d\sigma_A d\sigma_B)} \leq C^2\ \Vert h \Vert_{L^p(\R^{m+n})}.$$
\end{lemme}

This result states that one can deduce a restriction theorem on $A\times B$ from individual theorems on $A$ and $B$. This motivates the study of more general surfaces $S\subset \R^{m+n}$ for which, for all $y\in\R^{n}$, the level set $S_{y}:=\{x\in\R^{m} : (x,y) \in S\}$ admits a restriction theorem. Our idea is to apply this restriction inequality in the $x$ variable and then use dual Sobolev embeddings in $y$ to play the role of a second "restriction" inequality. This is the case of the cone $\mathcal{C}$, for which the level sets $\mathcal{C}_{\xi_{d+1}}$ are spheres of radius $\xi_{d+1}$, and for the surface $\mathcal{M}$, for which the level sets $\mathcal{M}_{\xi_{4}}$ are cones. We advise the reader to have a brief look at the proof of this lemma before looking at the proofs of Theorems 2, 3 and 4, since the structures of all the proofs are the same, but appear in a transparent way and without technicalities in the proof of Lemma 1.1.\\

What we consider to be the main asset of this approach is that it illustrates the role of curvature in restriction theorems. In the case of the cone, the way the spherical level sets are arranged along the vertical $\xi_{d+1}$-axis does not provide additional curvature to work with (since the spheres are linearly scaled), and thus applying the Tomas-Stein theorem on those spheres suffices to make use of all the available curvature. This is why this method gives a sharp result for the cone. For the surface $\M$ however, the conic level sets are disposed in a way that provides additional curvature along the $\xi_{4}$-axis. This additional curvature is of course ignored when applying a restriction theorem on a given level set. It can be taken into account by using a suitable Sobolev embedding, but then the Lebesgue norm in the restriction estimate becomes anisotropic and does not really provide a Strichartz estimate in its dual form, a phenomenon we detail in Remark 3.1. To keep the Lebesgue norm isotropic, we can only partially make use of the available curvature along the $\xi_{4}$-axis, and this is where the sharpness of Theorem 3 is lost.
\\

In Section 2 of this article we will prove Theorem 2 for the cone, before turning to the surface $\mathcal{M}$ and the proofs of Theorem 3 and Theorem 4 in Section 3, which will be concluded by Proposition 3.1. The proof of the previous lemma can be found in the Appendix. 

Throughout this note, we will denote by $d\sigma_R$ the surface measure on the sphere $\Spd_R\subset \R^d$ of radius $R$, defined by
$$ d\sigma_R(A) = \lambda(\{ta: a\in A, 0< t \leq 1\})$$
for any measurable subset $A \subset\Spd_R$, with $\lambda$ being the Lebesgue measure on $\R^d$. Finally, in the computations, $C$ will be a positive constant that may vary from line to line.

\section{Restriction theory on the cone: proof of Theorem 2}

In this section we will use the notations of Theorem 2. As mentioned in the introduction, the exponent $p_0$ is the endpoint one in Theorem 1, which applies to the sphere of $\R^d$. This is what motivates our technique, which consists in deducing Theorem 2 by applying Theorem 1 to the $(d-1)$-dimensional spheres of increasing radius that build up the $d$-dimensional cone. \\
We will thus need the following lemma, which is an extension of Theorem 1 to spheres of varying radius. Its proof is a simple scaling argument left to the reader.

\begin{lemme}
Let $\Spd_R$ be the sphere of radius $R>0$ included in $\R^d$, equipped with the surface measure $d\sigma_R$ (not normalized). There exists a constant $C>0$ (independent of R) such that for any function $f\in \s(\R^{d})$, one has:
$$ \Vert \hat{f}|_{\Spd_R} \Vert_{L^2(\Spd_R, d\sigma_R)} \leq C R^{\alpha} \Vert f \Vert_{L^{p_0}(\R^{d})},$$
with  $\alpha =  \frac{1}{p_0'} = \frac{d-1}{2(d+1)}$.
\end{lemme}

\begin{proof2}[Proof of Theorem 2.]

Consider $ f\in\s(\R^{d+1}) $. The first step to prove the theorem is to write $\Vert \hat{f}|_{\cont} \Vert_{L^2(\cont, d\sigma)}$ as an integral over spheres, using polar coordinates. We have:
\begin{align*}
    \Vert \hat{f}|_{\cont} \Vert^2_{L^2(\cont, d\sigma)}& = \itg_{\R^d} \big\vert\hat{f}(\xi,|\xi|) \big\vert^2 \frac{d\xi}{|\xi|}\\
    &= \itg_0^{+\infty} \itg_{\Spd_{\rho}} \vert \hat{f}(\xi,\rho) \vert^2 d\sigma_{\rho}(\xi) \frac{d\rho}{\rho}
\end{align*}
with $\Spd_{\rho}$ the sphere of radius $\rho$, and $d\sigma_{\rho}$ the measure on $\Spd_{\rho}$ as in Lemma 2.1.\\
Now, we want to use a restriction theorem on the sphere of radius $\rho$, that is to say apply Lemma 2.1. It yields

$$\itg_0^{+\infty} \itg_{\Spd_{\rho}} \vert \hat{f}(\xi,\rho) \vert^2 d\sigma_{\rho}(\xi) \frac{d\rho}{\rho}
    \leq C \itg_0^{+\infty} \bigg( \itg_{\R^d} \vert \Grave{f}(x,\rho) \vert^{p_0} dx \bigg)^{\frac{2}{p_0}} \frac{d\rho}{\rho^{1-2\alpha}}$$
with $\Grave{f}(x,\rho)$ the partial Fourier transform of $f$ in the last coordinate, that is:
$$\Grave{f} (x,\rho) = \int_{\R}e^{-ix_{d+1}\rho}\,f(x,x_{d+1})\,dx_{d+1}.$$
At this point, since $p_0 \leq 2$, we can use the Minkowski inequality: 
$$  \itg_0^{+\infty} \bigg( \itg_{\R^d} \vert \Grave{f}(x,\rho) \vert^{p_0} dx \bigg)^{\frac{2}{p_0}} \frac{d\rho}{\rho^{1-2\alpha}} 
\leq C \bigg(\itg_{\R^d} \bigg( \itg_0^{+\infty} \vert \Grave{f}(x,\rho) \vert^2 \frac{d\rho}{\rho^{1-2\alpha}} \bigg)^{\frac{p_0}{2}} dx \bigg)^{\frac{2}{p_0}} .$$
We recognize on the right-hand side a piece of the homogeneous Sobolev norm of $f(x,\cdot )$:
$$ \bigg( \itg_0^{+\infty} \vert \Grave{f}(x,\rho) \vert^2 \frac{d\rho}{\rho^{1-2\alpha}} \bigg)^{\frac{1}{2}} \leq \Vert f(x,\cdot )\Vert_{\dot H^s(\R)},$$
where
$$ \Vert f(x,\cdot )\Vert_{\dot H^s(\R)}=\bigg( \itg_{\R} \vert \Grave{f}(x,\rho) \vert^2 \vert \rho \vert^{2s}d\rho \bigg)^{\frac{1}{2}},$$
with $s = -\frac{1}{2} + \alpha$. Now, by the dual Sobolev embedding (that can be thought of as the second restriction-type inequality of Lemma 1.1), we find exactly $$\Vert f(x,\cdot )\Vert_{\dot H^s(\R)} \leq C \Vert f(x,\cdot ) \Vert_{L^{p_0}(\R)}. $$ 
Indeed, $s = - \frac{1}{d+1}$ and $\frac{1}{2} - \frac{1}{p_0} = \frac{d+1 - (d+3)}{2(d+1)} = s$ since $\alpha = \frac{1}{p_0'}$. Finally, 
\begin{align*}
    \bigg(\itg_{\R^d} \bigg( \itg_0^{+\infty} \vert \Grave{f}(x,\rho) \vert^2 \frac{d\rho}{\rho^{1-2\alpha}} \bigg)^{\frac{p_0}{2}} dx \bigg)^{\frac{2}{p_0}}
    &\leq C \bigg(\itg_{\R^d} \bigg( \itg_{\R} |f(x,y)|^{p_0} dy \bigg)^{\frac{p_0}{p_0}} dx \bigg)^{\frac{2}{p_0}} \\
    &\leq C \Vert f \Vert_{L^{p_0}(\R^{d+1})}^2,
\end{align*}
which concludes the proof.
\end{proof2}

\section{Linearized Euler equation in the rotational framework: proofs of Theorems 3 and 4}

We now turn our attention to the surface $\M$ and the proof of Theorem 3. The surface measure is now $d\sigma = \frac{d\xi}{|\xi|^2}$. We will reproduce the same technique we used in the previous section: indeed, notice that the level sets of $\M$ are cones of varying opening angles. Our goal is thus to write the $L^2$ norm of $f\vert_\M$ as an integral over all level sets and apply a suitable restriction theorem at each level.
We thus need to extend Theorem 2 into a corollary for cones of arbitrary angle:
\begin{cor}
Let $\cont_{\rho}:= \{ (\xi, \rho \vert \xi \vert),\ \xi\in\R^2 \}$ be a cone in  $\R^3$, with surface measure $d\sigma_\rho =\frac{d\xi}{\vert \xi \vert}$. Then, for any function $f\in \s(\R^{3})$, one has:
$$ \Vert \hat{f}|_{\cont_{\rho}} \Vert_{L^2(\cont_{\rho}, d\sigma_\rho)} \leq C \vert \rho \vert^{-\frac{1}{6}}\vert\  \Vert f \Vert_{L^{\frac 6 5}(\R^{3})}. $$
\end{cor}
This is directly deduced from the case $\rho=1$ by rescaling correctly $f$ in the third variable, with $\frac 6 5$ simply being the exponent in Theorem 2 for cones in $\R^3$.

We can now begin the proof of Theorem 3.

\begin{proof2}[Proof of Theorem 3.]
Consider $ f\in\s(\R^{4}) $ and set $\xi_h:=(\xi_1,\xi_2)$. We begin by writing $ \Vert \hat f\vert_{\M}  \Vert_{L^2(\M,d\sigma)} $ as an integral over all level sets by changing variables, letting $\mu=\dfrac{\xi_3}{\vert \xi\vert}$. By straightforward computations, we have the identities $d\xi_3=\vert \xi_h \vert (1-\mu^2)^{-\frac 3 2} d\mu $ and $\vert \xi \vert = \dfrac{\vert \xi_h\vert}{\sqrt{1-\mu^2}} $, which yield:
\begin{align*}
\Vert \hat f\vert_{\M}  \Vert^2_{L^2(\M,d\sigma)} &= \itg_{\R^3} \left\vert \hat f(\xi_h,\xi_3,\frac{\xi_3}{\vert \xi \vert } )\right\vert^2 \ \ \vert \xi\vert^{-2} d\xi_hd\xi_3 \\ &= \int_{-1}^{1}  \itg_{\R^2_{\xi_h}} \left\vert \hat f(\xi_h,  \dfrac{\mu}{\sqrt{1-\mu^2}} \vert \xi_h \vert , \mu )\right\vert^2 \ \frac{d\xi_h}{\vert \xi_h \vert} \ \dfrac{1}{\sqrt{1-\mu^2}} d\mu .
\end{align*}
We see that at a given level $\mu$, the integral over $\xi_h$ is simply a squared $L^2$ norm on a cone with the appropriate measure $\dfrac{d\xi_h}{\vert \xi_h \vert}$. We can apply Lemma 3.1 with $\rho=\dfrac{\mu}{\sqrt{1-\mu^2}}$ to get
\begin{align*}
\Vert \hat f\vert_{\M}  \Vert^2_{L^2(\M,d\sigma)} &\leq C\int_{-1}^{1}  \bigg( \itg_{\R^3_x} \left\vert  \Grave{f}(x, \mu )  \right\vert^{p} dx \bigg)^{\frac{2}{p}} \ \vert\mu\vert^{-\frac 1 3}(1-\mu^2)^{-\frac 1 3} d\mu ,
\end{align*}
with $p=\frac 6 5$, and $\Grave{f}(x,\mu)$ the Fourier transform of $f$ in the last coordinate $x_4$. As before, we now bound the $L^{2}_{\mu}L^{p}_{x}$ norm by a $L^{p}_{x}L^{2}_{\mu}$ norm using the Minkowski inequality. We are left with
\begin{align*}
\Vert \hat f\vert_{\M}  \Vert^2_{L^2(\M,d\sigma)} &\leq C  \bigg(\itg_{\R^3_x}\bigg( \int_{-1}^{1} \left\vert  \Grave{f}(x, \mu )  \right\vert^{2} \vert\mu\vert^{-\frac 1 3}(1-\mu^2)^{-\frac 1 3} d\mu \bigg)^{\frac{p}{2}} \  dx\bigg)^{\frac{2}{p}}.
\end{align*}
We now want to use Sobolev embeddings for the inner integral to bound it by a $L^{p}$ norm. Note that we have three singularities in $\mu=-1$, $0$, $1$, however near each of those points the term $\vert\mu\vert^{-\frac 1 3}(1-\mu^2)^{-\frac 1 3}=\vert\mu\vert^{-\frac 1 3}(1-\mu)^{-\frac 1 3}(1+\mu)^{-\frac 1 3}$ behaves like a simple Sobolev weight. We make use of this by splitting the integral in three:
$$\int_{-1}^{1} \left\vert  \Grave{f}(x, \mu )  \right\vert^{2} \vert\mu\vert^{-\frac 1 3}(1-\mu^2)^{-\frac 1 3} d\mu = \int_{-1}^{-\frac 1 2} ... + \int_{-\frac 1 2}^{\frac 1 2} ... + \int_{\frac 1 2}^{1} ...\ .$$
Let us bound the last term, the others being dealt with identically: we simply use that, for $\frac 1 2 \leq \mu \leq 1$, we have $\vert\mu\vert^{-\frac 1 3}(1-\mu)^{-\frac 1 3}(1+\mu)^{-\frac 1 3}\leq C (1-\mu)^{-\frac 1 3}$. This yields, making an appropriate change of variables:
\begin{align*}
\int_{\frac 1 2}^{1}\left\vert  \Grave{f}(x, \mu )  \right\vert^{2} \vert\mu\vert^{-\frac 1 3}(1-\mu^2)^{-\frac 1 3} d\mu &\leq C\int_{\frac 1 2}^{1}\left\vert  \Grave{f}(x, \mu )  \right\vert^{2}  (1-\mu)^{-\frac 1 3} d\mu \\ &=C\int_{0}^{\frac 1 2}\left\vert  \Grave{f}(x, 1-\nu )  \right\vert^{2}  \nu^{-\frac 1 3} d\nu .
\end{align*}
The Sobolev exponent we seek is $s=-\frac{1}{3}$, to exploit that $L^{p}$ is embedded in $\dot{H}^{s}(\R)$. Since $\vert \nu \vert \leq 1$, we can bound $\nu^{-\frac 1 3}$ by $\nu^{-\frac 2 3}$. Note that this is very blunt, and most of the sharpness of the result is lost due to this (we further discuss this in Remark 3.1). Using this bound anyway, we get
\begin{align*}
\int_{\frac 1 2}^{1}\left\vert  \Grave{f}(x, \mu )  \right\vert^{2} \vert\mu\vert^{-\frac 1 3}(1-\mu^2)^{-\frac 1 3} d\mu &\leq C\int_{0}^{\frac 1 2}\left\vert  \Grave{f}(x, 1-\nu )  \right\vert^{2}  \nu^{-\frac 2 3} d\nu \\&\leq C\int_{\R_{\nu}}\left\vert  \Grave{f}(x, 1-\nu )  \right\vert^{2}  \vert\nu\vert^{-\frac 2 3} d\nu\\&\leq C \Vert f(x,\cdot )\Vert ^{2}_{L^{p}(\R)}.
\end{align*}
In the last step we used Sobolev embeddings, as well as the fact that the inverse Fourier transform in $\nu$ of $\Grave{f}(x, 1-\nu )$ is $f$ up to a phase and a sign in the argument, both making no difference in the $L^p$ norm.
The same bound holds for the two other terms, so that the whole integral over $\mu$ is controlled by $\Vert f(x,.)\Vert ^{2}_{L^{p}(\R)}$, which in turns yields
\begin{align*}
\Vert \hat f\vert_{\M}  \Vert^2_{L^2(\M,d\sigma)} &\leq C  \bigg(\itg_{\R^3_x}\bigg(  \Vert f(x,.)\Vert ^{2}_{L^{p}(\R)}\bigg)^{\frac{p}{2}} \  dx\bigg)^{\frac{2}{p}}\\&= C  \Vert f\Vert ^{2}_{L^{p}(\R^4)}.
\end{align*}
This concludes the proof of Theorem 3.
\end{proof2}

\begin{rmq}
Let us briefly discuss the sharpness of this result.
As mentioned in the proof, we used at some point the blunt bound $\nu^{-\frac{1}{3}} \leq \nu^{-\frac{2}{3}}$, because the Sobolev embedding $L^{\frac 6 5}(\R) \hookrightarrow \dot H^{-\frac 1 3}(\R)$ naturally provides a bound by an $L^{\frac 6 5}$ norm in the last coordinate, so that we obtain in the end a $L^{\frac 6 5}$ norm on the whole space $\R^4$. However, this obviously makes us lose precision on the result. To avoid this, we could have kept the term $\nu^{-\frac 1 3}$, and instead use the Sobolev embedding $L^{\frac 3 2}(\R) \hookrightarrow \dot H^{-\frac 1 6}(\R)$. However, at the end, we would have got the following restriction result involving anisotropic Lebesgue norms:
$$ \Vert \hat f\vert_{\M}  \Vert_{L^2(\M,d\sigma)} \leq C \Vert f \Vert_{L^{\frac{6}{5}}(\R^{3}_{x_1,x_2,x_3},L^{\frac{3}{2}}(\R_{x_4}))}. $$
In terms of estimates for the semigroup $e^{it\frac{D_3}{\vert D \vert}}$, it writes
$$\Vert e^{it\frac{D_3}{\vert D \vert}} u_0 \Vert_{L^6_{x}L^3_{t}} \leq C \Vert u_0 \Vert_{\dot H^1(\R^3)}.$$
This is not really a Strichartz estimate, we would have liked a $L^3_t L^6_x$ norm instead, but to our knowledge there is no way around this issue (see \cite{bbg} where the same phenomenon occurs). We know thanks to \cite{cor} that $$ \Vert e^{it\frac{D_3}{\vert D \vert}} u_0 \Vert_{L^3_{t}L^6_{x}} \leq C \Vert u_0 \Vert_{\dot H^1(\R^3)} $$ is true, and is sharp in the sense that for a $L^6_x$ norm in space, we cannot do better than a $L^3_t $ norm in time.
\\The source of the anisotropy in the result is due to the curvature of $\M$ along the $\xi_4$ coordinate, which is completely ignored when applying Theorem 2 on each conic level set. This additional curvature only plays a role in the integration over $\xi_4$, and allows us to sharpen the Sobolev embedding we use, but at the cost of the Lebesgue exponents not matching in the end. \\
A way to try and obtain a sharp Strichartz estimate would be to aim for the same exponents in time and space (so that their order does not matter). This would be the following estimate $$ \Vert e^{it\frac{D_3}{\vert D \vert}} u_0 \Vert_{L^4_{t}L^4_{x}} \leq C \Vert u_0 \Vert_{\dot H^{\frac 3 4}(\R^3)} $$ from \cite{cor}. Those exponents are too strong to be achieved using cone restriction theorems applied to the level sets, since this would require Theorem 2 but with $p_0=\frac 4 3$ in dimension $d=2$ which is not possible, see \cite{tao}.
\end{rmq}

We shall now prove Theorem 4, we will give less details since the proof is very similar to the previous one.

\begin{proof2}[Proof of Theorem 4.]
Let $f\in\s(\R^4)$, $F \in (1,\infty)$ and $\xi_h = (\xi_1,\xi_2)$. Recall that $\vert \xi \vert_{F}^{2}=\vert\xi_{h}\vert^{2}+F^{2}\vert\xi_{3}\vert^{2}$ and $\mathcal{M}_{F}=\left\{ \bigg(\xi,\dfrac{\vert \xi \vert_{F}}{F\vert\xi\vert}\bigg), \xi \in \R^{3}\right\}\subset \R^4,$ thus we want to compute
$$ \itg_{\R^3} \left\vert \hat f \bigg(\xi, \dfrac{\vert \xi \vert_{F}}{F\vert\xi\vert}\bigg) \right\vert^2 \frac{d\xi}{\vert \xi\vert^2}\ \cdotp$$
For a bit more simplicity, we shall restrict our attention to half of this surface, say $\xi_3>0$. It does not change the result since the surface is symmetric relatively to the axis $\xi_3=0$.\\

The first step is to write $\mu =  \dfrac{\vert \xi \vert_{F}}{F\vert\xi\vert}$, thus we can express $\xi_3$ in terms of $\xi_h$ and $\mu$. It yields 
$$ \xi_3 = \sqrt{ \frac{\mu^2-F^{-2}}{1-\mu^2} }\ \ \vert\xi_h\vert,\ \ \ \ \ \ \ \ \ \ \ \ \ \ \ \ \ \ \ \ \ d\xi_3 = \mu \frac{1-F^{-2}}{(1-\mu^2)^{\frac 3 2}(\mu^2-F^{-2})^{\frac 1 2 }}\ \ \vert\xi_h\vert\ \ d\mu\ \ \ \ \ \ \ \ \ \ \ \ \ \ \ \ \text{and}\ \ \ \ \ \ \ \ \ \ \ \ \ \ \ \ \ \ \vert\xi\vert^2 = \vert\xi_h\vert^2\bigg(\frac{1-F^{-2}}{1-\mu^2}\bigg).$$
Then, straightforward computations give 
$$ \itg_{\R^2_{\xi_h}\times\R^+_{\xi_3}} \left\vert \hat f (\xi, \dfrac{\vert \xi \vert_{F}}{F\vert\xi\vert}) \right\vert^2 \frac{d\xi}{\vert \xi\vert^2}
= \itg_{\frac{1}{F}}^1 \itg_{\R^2_{\xi_h}} \left\vert \hat f \bigg(\xi_h, \sqrt{ \frac{\mu^2-F^{-2}}{1-\mu^2} }\ \ \vert\xi_h\vert , \mu \bigg) \right\vert^2  \frac{\mu}{(1-\mu^2)^{\frac 1 2}(\mu^2-F^{-2})^{\frac 1 2 }}\ \ \frac{d\xi_h}{\vert \xi_h\vert} d\mu.  $$
Now, we are ready to apply Corollary 3.1 with $\rho = \sqrt{ \frac{\mu^2-F^{-2}}{1-\mu^2} }$ and $p= \frac{ 6}{5}$:
\begin{align*}
    \itg_{\frac{1}{F}}^1 \itg_{\R^2_{\xi_h}} 
    \left\vert \hat f \bigg(\xi_h, \sqrt{ \frac{\mu^2-F^{-2}}{1-\mu^2} }\ \ \vert\xi_h\vert , \mu \bigg) \right\vert^2 & \frac{\mu}{(1-\mu^2)^{\frac 1 2}(\mu^2-F^{-2})^{\frac 1 2 }}\ \ \frac{d\xi_h}{\vert \xi_h\vert} d\mu\\
    &\leq C \itg_{\frac{1}{F}}^1 \bigg(\itg_{\R^3_x} \vert \Grave{f}(x,\mu) \vert^p dx \bigg)^{\frac 2 p} \mu (1-\mu^2)^{-\frac 1 3}(\mu^2-F^{-2})^{-\frac 2 3}d\mu.
\end{align*}
The Minkowski inequality turns the $L^2_{\mu}L^p_x$ norm into a $L^p_x L^2_{\mu}$ norm:

\begin{align*}
     \itg_{\frac{1}{F}}^1 \bigg(\itg_{\R^3_x} \vert \Grave{f}(x,\mu) \vert^p dx \bigg)^{\frac 2 p}& \mu (1-\mu^2)^{-\frac 1 3}(\mu^2-F^{-2})^{-\frac 2 3}d\mu\\ 
     &\leq \Bigg( \itg_{\R^3_x} \bigg( \itg_{\frac{1}{F}}^1 \vert \Grave{f}(x,\mu) \vert^2 \mu (1-\mu^2)^{-\frac 1 3}(\mu^2-F^{-2})^{-\frac 2 3}d\mu \bigg)^{\frac p 2 } dx \Bigg)^{\frac 2 p}.
\end{align*}

As before, it is enough to prove that $$\itg_{\frac{1}{F}}^1 \vert \Grave{f}(x,\mu) \vert^2 \mu (1-\mu^2)^{-\frac 1 3}(\mu^2-F^{-2})^{-\frac 2 3}d\mu \leq C_F \Vert f(x,y) \Vert_{L^p(\R_y)}^2.$$
To do so, we shall split the previous integral around the singularities at $F^{-1}$ and $1$, and then use the Sobolev embedding $L^{\frac 6 5}(\R) \hookrightarrow \dot H^{-\frac 1 3}(\R)$. Indeed,
$$\itg_{\frac{1}{F}}^1 \vert \Grave{f}(x,\mu) \vert^2 \mu (1-\mu^2)^{-\frac 1 3}(\mu^2-F^{-2})^{-\frac 2 3}d\mu  = \itg_{\frac{1}{F}}^{ \frac{F+1}{2F} } ... + \itg_{ \frac{F+1}{2F} }^1 ... \ .$$
For $\mu \in [F^{-1}, \frac{F+1}{2F}]$, we have $\mu (1-\mu^2)^{-\frac 1 3}(\mu^2-F^{-2})^{-\frac 2 3} \leq \frac{F}{(F-1)^{\frac 1 3}} (\mu - F^{-1})^{\frac 2 3}$, and this is exactly the Sobolev weight we wanted (the way we deal with this term is exactly the same than for Theorem 3). \\
For $\mu \in [\frac{F+1}{2F},1]$, we shall use the inequality $(\mu-F^{-1})^{-\frac 2 3} \leq \bigg(\frac{F}{F-1}\bigg)^{\frac 1 3} (1-\mu)^{-\frac 1 3}$. It yields
$$\mu (1-\mu^2)^{-\frac 1 3}(\mu^2-F^{-2})^{-\frac 2 3} \leq C \bigg(\frac{F}{F-1}\bigg)^{\frac 1 3} (1-\mu)^{-\frac 2 3},$$
once again it is the desired Sobolev weight. It gives
$$ \itg_{\frac{1}{F}}^1 \vert \Grave{f}(x,\mu) \vert^2 \mu (1-\mu^2)^{-\frac 1 3}(\mu^2-F^{-2})^{-\frac 2 3}d\mu \leq C \frac{F}{(F-1)^{\frac 1 3}}\bigg(\itg_{\R_y} \vert f(x,y)\vert^p dy \bigg)^{\frac 2 p},$$
which concludes the proof.
\end{proof2}

\begin{rmq}
The fact that the estimates blow up as $F\rightarrow 1$ is very intuitive : the surface then becomes the hyperplane $\{\xi_4 = 1\}$. We refer to \cite{fun} for a study of the primitive equations in the case $F=1$ with viscosity. However the blow-up as $F\rightarrow\infty$ is a bit less intuitive, since the estimate still holds for $F= \infty$. We believe this is due to the fact that $\M_F$ is located above the hyperplane ${\xi_{4}=\frac{1}{F}}$: for finite $F$, the singularity at $0$ of the Sobolev weight in the surface measure does not intervene, whereas it appears in the case $F=\infty$ and offers a gain of integrability at the origin.\\
Finally, the Sobolev weight that appears in the Jacobian is exactly the one wanted at one of the two singularities. Thus, contrary to the surface $\M$, in this case we could not use a sharper Sobolev embedding to improve the estimate, even if we allowed anisotropic norms. This might be a sign that the result is in some way sharper than for Theorem 3, but we did not investigate this further.
\end{rmq}

To conclude, we shall say a few words about the need of Sobolev norms instead of $L^2$ norms. As said in the introduction, the surface studied is very flat at infinity, and we have seen that the hyperplane can never admit non-trivial restriction theorems. The following proposition is the consequence of the flatness and the non-compactness of $\M$, it is stated in dimension $d$ for a bit more generality.

\begin{prop}
Consider the $d$-dimensional surface $$ \M:= \left\{ \bigg(\xi, \frac{\xi_d}{\vert \xi \vert}\bigg),\ \ \ \xi \in \R^d \right\},$$
endowed with the measure $d\sigma(\xi,\frac{\xi_d}{\vert \xi\vert}) = d\xi$.  Let $p>1$. Then one cannot find an exponent $q$ and a constant $C$ such that for any $f\in \s(\R^{d+1})$,
$$ \Vert \hat f \vert_{\M} \Vert_{L^q(\M,d\sigma)} \leq C \Vert f \Vert_{L^p(\R^{d+1})}.$$
\end{prop}

\begin{proof2}[Proof of Proposition 3.1.]
Consider  $f(x\dpu) = \frac{1}{\sqrt{1 + x\dpu^2}}$, and $\psi \in \s(\R^d)$ such that $$\supp \hat \psi \subset \T: = \{ (\xi',\xi_d)\in \R^{d-1} \times \R,\ \ \ \vert \xi' \vert \leq 1,\ \ \ \vert \xi_d \vert \leq  \vert \xi' \vert \}.$$
Let $h(x_1,...,x\dpu) = \psi(x_1,...,x_d) f(x\dpu)$, and let $$h_R =  \psi(Rx_1,..., Rx\dmu, \frac{x_d}{R^{d-1}}) f(x\dpu).$$
\\
Then,  
$$\widehat{h_R} =  \hat \psi(\frac{\xi_1}{R},..., \frac{\xi\dmu}{R},{R^{d-1}}{\xi_d}) \hat f(\xi\dpu).$$
This can be written, with obvious notation, as $\widehat{h_R} = \widehat{\psi_R} \hat f$. Moreover, $\widehat{\psi_R}$ is supported in $$ \T_R: = \{ (\xi',\xi_d)\in \R^{d-1} \times \R,\ \ \ \vert \xi' \vert \leq R,\ \ \ \vert \xi_d \vert \leq  \frac{\vert \xi' \vert}{R^{d-1}} \}.$$ 
Then we shall study the norms of $h_R$ and $\widehat{h_R}\vert_{\M}$. \\
On the one hand, $ \Vert h_R \Vert_{L^p(\R^{d+1})} = \Vert h \Vert_{L^p(\R^{d+1})} $ whenever $p>1$.\\
On the other hand, using the fact that $\frac{\xi_d}{\vert\xi\vert}$ is very close to 0 for $\xi\in\T_R$, we find 
$$  \Vert \widehat{h_R}\vert_{\M} \Vert_{L^q(\M)} \geq  \inf\limits_{\tau \in [-\frac{1}{R^{d-1}}, \frac{1}{R^{d-1}}]} \vert \hat f (\tau) \vert \ \ \Vert \widehat{\psi} \Vert_{L^q(\R^d)} .$$
Yet the infimum of $\hat f$ blows up as $R \rightarrow \infty$, hence one cannot have an inequality such as $$  \Vert \hat u \vert_{\M} \Vert_{L^q(\M)} \leq C \Vert u \Vert_{L^p(\R^{d+1})}$$
holding for all $u\in\s(\R^{d+1})$.
\end{proof2}

\section*{Appendix}

\begin{proof2}[Proof of Lemma 1.1.]
Let $f\in \s(\R^{m+n})$. Applying the restriction inequality on $B$, we have:
\begin{align*}
    \Vert \hat{f}|_{A\times B} \Vert^q_{L^q(d\sigma_A d\sigma_B)} &= \itg_A \itg_B |\hat{f}(\xi,\tau)|^q  d\sigma_B(\tau)d\sigma_A(\xi) \\
    &\leq C^q \itg_A \bigg( \itg_{\R^n} \vert \Acute{f}(\xi,y) \vert^p dy \bigg)^{\frac{q}{p}} d\sigma_A(\xi),
\end{align*}    
with $\Acute{f}(\xi,y)$ the partial Fourier transform of $f$ in the $m$ first coordinates. Now, by the Minkowski inequality, since $p\leq q$,
$$\itg_A \bigg( \itg_{\R^n} \vert \Acute{f}(\xi,y) \vert^p dy \bigg)^{\frac{q}{p}} d\sigma_A(\xi) \leq   \bigg( \itg_{\R^n} \bigg( \itg_A \vert \Acute{f}(\xi,y) \vert^q  d\sigma_A(\xi)  \bigg)^{\frac{p}{q}} dy\bigg)^{\frac{q}{p}}.$$

Finally, it suffices to apply the restriction inequality on $A$:
\begin{align*}
     C^q \bigg( \itg_{\R^n} \bigg( \itg_A \vert \Acute{f}(\xi,y) \vert^q  d\sigma_A(\xi)  \bigg)^{\frac{p}{q}} dy\bigg)^{\frac{q}{p}} 
   & \leq C^{2q} \bigg( \itg_{\R^n} \bigg( \itg_{\R^m} \vert {f}(x,y) \vert^p  dx \bigg)dy \bigg)^{\frac{q}{p}} \\
    & \leq C^{2q} \Vert f \Vert^q_{L^p(\R^{m+n})},
\end{align*}
which concludes the proof.
\end{proof2}

\end{document}